\documentclass[11pt]{amsart}

\usepackage {amsfonts,amsthm}
\usepackage[english]{babel}

\theoremstyle{plain}
\newtheorem{Th}{Theorem}
\newtheorem*{Pro}{Proposition}
\newtheorem*{Cor}{Corollary}
\newtheorem{Def}{Definition}

\def\card{{\rm{card}}}
\def\A{almost periodic multiple discrete set}
\def\C{{\mathbb{C}}}
\def\R{{\mathbb{R}}}
\def\Z{{\mathbb{Z}}}
\def\N{{\mathbb{N}}}

\def\s{\sigma}
\def\g{\gamma}
\def\d{\delta}
\def\f{\varphi}
\def\a{\alpha}
\def\l{\lambda}
\def\b{\beta}

\def\e{{\varepsilon}}
\def\G{\Gamma}
\def\supp{{\rm{supp}}}
\def\mes{{\rm{mes}}}

\def\ws{~\hfill$\Box$}

\def\be{\begin{equation}}
\def\ee{\end{equation}}

\begin{document}
\title[Perturbations of discrete lattices]
{Perturbations of discrete lattices \\ and almost periodic sets}

\author[Favorov S., Kolbasina Ye.]{Favorov S., Kolbasina Ye.}

\address{Mathematical School, Kharkov National University, Swobody sq.4,
Kharkov, 61077 Ukraine}

 \email{Sergey.Ju.Favorov@univer.kharkov.ua,kvr${\_}$jenya@mail.ru}

\date{}

\begin{abstract}
 A discrete set in the $p$-dimensional Euclidian space is {\it almost periodic}, if  the
 measure with the unite masses at points of the set is almost periodic in the weak sense.
We propose to construct positive almost periodic discrete sets as an almost periodic
perturbation of a full rank discrete lattice.  Also we prove that each almost periodic
discrete set on the real axes is an almost periodic perturbation of some arithmetic
progression.

Next,  we consider signed almost  periodic discrete sets, i.e., when the signed
measure with masses $\pm1$ at points of a discrete set is almost periodic. We
construct a signed discrete set that is not almost periodic, while the corresponding
signed measure is almost periodic in the sense of distributions. Also, we construct a
signed almost periodic discrete set such that the measure with masses $+1$ at all
points of the set is not almost periodic.
\end{abstract}



\maketitle

The concept of  almost periodicity  plays an important role in various branches of
 analysis. In  particular, almost periodic discrete sets are used for investigation of zero sets
of some holomorphic functions (cf.~\cite{L},\cite{FRR}), in  value distribution theory
of some classes of meromorphic functions (cf.~\cite{F}), as a model of quasicrystals
(cf.~\cite{La},\cite{M}). Note that in \cite{La} the question (Problem 4.4) was raised
if there exist other discrete almost periodic sets in $\R^p$, besides of the form
$L+E$ with a discrete lattice $L$ and a finite set $E$.

Here we propose several ways to construct almost periodic discrete sets. Next, we
introduce signed almost periodic discrete sets. To formulate our result beforehand we
have to recall some known definitions (see, for example, \cite{C}).

\medskip
A continuous function $f(x)$ in $\R^p$ is {\it almost periodic}, if for any $\e>0$ the
set of $\e$-almost periods of $f$
 $$
  \{\tau\in\R^p:\,\sup_{x\in\R^p}|f(x+\tau)-f(x)|<\e\}
  $$
  is a relatively dense set in $\R^p$. The latter means that there is
  $R=R(\e)<\infty$ such that any ball of radius $R$ contains an $\e$-almost period of
  $f$.

Note that  almost periodic functions are uniformly bounded in $\R^p$. Besides, every
almost periodic function is the uniform in $x\in\R^p$ limit of a sequence of
exponential polynomials of the form
 \be\label{pol}
 P(x)=\sum_m c_m e^{i\langle x,\l_m\rangle},\quad \l_m\in\R^p,\quad  c_m\in\C,
  \ee
 here $\langle .,.\rangle$ is the scalar product in $\R^p$.

\medskip
A Borel measure $\mu$ in $\R^p$ is {\it almost periodic} if it is almost periodic in
the weak sense, i.e., for any continuous function $\f$ in $\R^p$ with a compact
support the convolution \be
    \label{conv}
    \int\f(x+t)~d\mu(t)
\ee
  is an almost periodic function in $x\in\R^p$ (see \cite{R}).

The definition suits to signed measures as well.

\begin{Th}[\cite{R}, Theorems 2.1 and 2.7]\label{R}
For any signed almost periodic measure $\mu$ in $\R^p$ there exists $M<\infty$ such
that the variation $|\mu|$ satisfies the condition
\begin{equation}
       \label{b1}
         |\mu|(B(c,1))< M \quad \forall c\in\R^p.
\end{equation}
Besides, there exists uniformly in $x\in\R^p$ a finite limit
 $$
        D(\mu)= \lim_{R\to\infty}\frac{\mu(B(x,R))}{\omega_pR^p}.
 $$
\end{Th}
Here $B(x,R)$ is an open ball with the center at the point $x$ and radius $R$,
$\omega_p$ is the volume of $B(0,1)$.

\medskip
We will consider some generalization of discrete sets, namely multiple discrete sets
in $\R^p$. This means that a number $m(x)\in\N$ corresponds to each point $x$ from a
discrete set. We denote this object by $A=\{(x,m(x))\}$ and the corresponding discrete
set by $s(A)$. Also, we will write $A=(a_k)$, where every term $a\in s(A)$ appears
$m(a)$ times in the sequence $(a_k)$.

In the case $p=2$ the definition coincides with the definition of the  divisor of an
entire function in the complex plane.

\medskip
The definition of almost periodic sets has an evident generalization to \A s.
\begin{Def}\label{m} A  multiple discrete set $A$ is {\it almost periodic}, if its associate
measure
 \be\label{f}
 \mu_A=\sum_{x\in s(A)} m(x)\d_x,
 \ee
 where $\d_x$ is the unit mass at the point $x$, is almost periodic.
\end{Def}
 Put
  $$
  \card(A\cap E)=\sum_{x\in s(A)\cap E}m(x)
  $$
 for any $E\subset\R^p$. The following result is a consequence of Theorem \ref{R}.

 \begin{Th}\label{RF}
For any \A\,$A$ there exists $M<\infty$ such that
\begin{equation}
       \label{b2}
         \card\;(A\cap B(x,1))< M \quad \forall  x\in\R^p.
\end{equation}
Besides, there exists uniformly in $x\in\R^p$ a finite density
\begin{equation} \label{denc}
        D(A)= \lim_{R\to\infty}\frac{\card(A\cap B(x,R))}{\omega_pR^p}.
\end{equation}
\end{Th}
Another proof of Theorem \ref{RF} see in \cite{FK}.

\medskip
 There is a geometric criterium for multiple discrete sets to be
almost periodic.
\begin{Th}[\cite{FK}, Theorem 11]\label{FK}
 An \A\  $(a_n)\subset\R^p$ is almost periodic if and only if for each $\e>0$
the set of $\e$-almost periods of $(a_n)$
 \be\label{bij}
\{\tau\in\R^p:\,\exists \hbox{\ a bijection}\, \s:\,\N\to\N\quad\hbox{ such that
}\quad |a_n+\tau-a_{\sigma(n)}|<\e\quad \forall\,n\in\N\}
 \ee
 is relatively dense in $\R^p$.
\end{Th}

At the first time  almost periodic divisors appeared in papers \cite{L} and \cite{T},
where only shifts along real axis were considered. The definition of almost
periodicity based on the above geometric property. An analog of Theorem \ref{FK} was
proved in \cite{FRR}.

\bigskip {\bf Almost periodic perturbations of discrete lattices.} Let
$F(x)=(F_1(x),\dots,F_p(x))$ be a mapping from $\R^p$ to $\R^p$ with almost periodic
components $F_j(x)$. For convenience of a reader, prove the following known assertion.

\begin{Pro}
For any $\e>0$ there exists a relatively dense set of common $\e$-almost periods of
 $F_j$ with integer coordinates.
\end{Pro}
{\bf Proof.} By the well-known Kronecker Theorem, the system of inequalities
 \be\label{s}
|\exp\langle\tau,\l_n\rangle-1|<\d,\quad \quad n=1,\dots,N,
 \ee
has a relatively dense in $\R^p$ set of solutions $\tau$ for any $\d>0$ and any
$\l_n\in\R^p,\;n=1,\dots,N$. Let $\{e_j\}_{j=1}^p$ be the natural basis in $R^p$.
Common solutions of (\ref{s}) and the system
 \be\label{s1}
|\exp2\pi\langle\tau,e_j\rangle-1|<\eta,\quad j=1,\dots,p,
 \ee
form a relatively dense set too. Whenever $\tau$ satisfies (\ref{s1}), there is
$r\in\Z^p$ such that $|\tau-r|<p\eta/\pi$. Hence for sufficiently small $\eta$ there
exists a relatively dense set of solutions $r\in\Z^p$ of system (\ref{s}) with $2\d$
instead of $\d$.
 If each function $F_j$ is an exponential polynomial of form (\ref{pol}),
 then there exist $\d$ and $\l_1,\dots,\l_N$ such that these solutions are common $\e$-almost
 periods of $F_j$. In the general case we can approximate the functions $F_j(x)$ by
 sequences of exponential polynomials. \ws

\medskip
Let $L$ be an arbitrary discrete full rank lattice in $\R^p$. Rewrite it in the form
$L=\{k\G,\,k\in\Z^p\}$, where $\G$ is a non-generated $p\times p$ matrix.  If
$r\in\Z^p$ is a common $\e$-almost period of components of $F$, then $\tau=r\G$ is an
$p\e$-almost period of the discrete multiple set
 \be\label{rep}
 A=(a_k),\quad a_k=k\G+F(k),\quad k\in\Z^p,
 \ee
 where the bijection $\s:\,Z^p\to\Z^p$ in (\ref{bij})
has the form $\s(k)=k+r$. Whenever all components $F$ are almost periodic, Theorem
\ref{FK} implies that $A$ is an almost periodic multiple discrete set. Also, note that
in the case of sufficiently small $\sup_{\R^p}|F(x)|$ we obtain an almost periodic
Delone set.

\medskip
It is easy to construct an almost periodic set in $\R^p$ without any periods. Take
$F(x)=(1/5)(\sin x_1,\dots,\sin x_p)$ for $x=(x_1,\dots,x_p)$, and put
$A=k+F(k),\,k\in\Z^p$. If $\tau\in\R^p\setminus\{0\}$ is a period of $A$, then
$k+\tau+F(k)=k'+F(k')$ for all $k\in\Z^p$ and some $k'=k'(k,\tau)$. Clearly,
$\tau=\tau^{(1)}+\tau^{(2)}$, where
$\tau^{(1)}=(\tau^{(1)}_1,\dots,\tau^{(1)}_p)\in\Z^p$ and modula of all components of
$\tau^{(2)}=(\tau^{(2)}_1,\dots,\tau^{(2)}_p)$ is less then $1/2$. Therefore,
$k'=k+\tau^{(1)}$ and $\sin(k_j+\tau^{(1)}_j)=\tau^{(2)}_j+\sin k_j$ for all
$k=(k_1,\dots,k_p)\in\Z^p$ and $j=1,\dots,p$, that is impossible.

Clearly, all vectors from a lattice $L$ are periods of every almost periodic set of
the form $L+E$ with a finite set $E$. Therefore we obtain an answer on the question
raised in \cite{La}.

\medskip
  Note  that in \cite{DF} and \cite{Ko} we get representation (\ref{rep})
 with a squire lattice $L$ and a {\it bounded} mapping $F(x)$
for a wide class of multiple discrete sets in $\R^p$, in particular, for every
 \A s. I do not know if each almost
periodic set has representation (\ref{rep})  with some lattice $L$ and a mapping
$F(x)$ with  almost periodic coordinates $F_j$. But this is true for almost periodic
sets in the real axis.

\begin{Th}\label{3}
Let $A=(a_k)_{k\in\Z}$, where $a_k\le a_{k+1}$ for all $k$, be an almost periodic
multiple discrete set in $\R$ with the density $D$. Then $a_k=Dk+f(k)$ with an almost
periodic function $f$.
 \end{Th}
In \cite{L} a similar result was obtained for real parts of zeros of almost periodic
entire functions from some special class.
\medskip

{\bf Proof.} Without loss of generality we may suppose that
 density $D$ of the set $A$ is equal to $1$. Also we suppose that $0\in A$ and $a_0=0$.
Take arbitrary $x,\,y,\,h\in\R,\,x<y,\,h>0$. Using Theorem \ref{FK}, take $L>2$ such
that any interval $i\subset\R$ of length $L$ contains a $1$-almost period $\kappa$ of
$A$. Since $\overline{i}-\kappa\subset (-L,\,L)$, we get
 $$
 \card(A\cap\overline{i})\le M,\quad\hbox{where}\quad M=\card(A\cap(-L-1,\,L+1)).
 $$
Take a $1$-almost period $\kappa$ of the set $A$
 such that $y<x+\kappa<y+L$. By definition, there is a bijection $\rho$ between all
points of the set $A\cap(x+\kappa,\,x+\kappa+h]$ and some points of the set
$A\cap(x-1,\,x+h+1]$. Moreover, the same $\rho$ is a bijection between some points of
the set $A\cap(x+\kappa,\,x+\kappa+h]$ and all points of the set $A\cap(x+1,\,x+h-1]$.
Therefore we have
 $$
|\card(A\cap(x+\kappa,\,x+\kappa+h])-\card(A\cap(x,\,x+h])|\le\card(A\cap((x-1,\,x+1]
\cup(x+h-1,\,x+h+1])\le2M.
  $$
 Since
   $$
  (x+\kappa,x+\kappa+h]\setminus(y,y+h]\subset(y+h,y+h+L),
  \quad(y,y+h]\setminus(x+\kappa,x+\kappa+h]\subset(y,y+L),
  $$
  we obtain
   $$
 |\card(A\cap(x,\,x+h])-\card(A\cap(y,\,y+h])|\le2M+\card(A\cap[(y,y+L)\cup(y+h,y+h+L)])\le 4M.
   $$
For any $T\in\N$ the half--interval $(x,\,x+Th]$ is the union of half--intervals
$(x+(j-1)h,\,x+jh],\,j=1,\dots,T$. If we set $y=x+(j-1)h,\,j=2,\dots,T$,
  we get
  $$
  |\card(A\cap(x,\,x+h])-T^{-1}\card(A\cap(x,\,x+Th])|\le 4M.
  $$
Taking into account (\ref{denc}) with $D(A)=1$, we obtain
 \be\label{est}
|\card(A\cap(x,\,x+h])-h|\le 4M\qquad\forall\,x\in\R,\,h>0.
 \ee
Furthermore, put
 $$
 n(t)=\card(A\cap(0,\,t])\quad\hbox{for}\quad
t>0,\quad n(t)=-\card(A\cap(t,\,0])\quad\hbox{for}\quad t<0,\quad n(0)=0.
 $$
Take $\e<L/(24M)$. Let $\tau>2L$ be an $\e$-almost period of the set $A$. Clearly, if
$x,\,y,\,x+\tau,\,y+\tau$ do not belong to  the $\e$-neighborhood $U_\e$ of the set
$s(A)$, then we have
 $$
n(y+\tau)-n(x+\tau)=n(y)-n(x).
 $$
 Therefore, the function $n(x+\tau)-n(x)$ takes the same number $p\in\N$ for
all $x\in\R\setminus(U_\e\cup(U_\e-\tau))$.

Next, denote by $E[a]$ the integer part of a real number $a$. Put $N=E[L/(4M\e)]+1$.
It is easily shown that
  \be\label{est-N}
 \frac{L}{4M\e}<N<\frac{\tau}{2M\e(E[\tau/L]+1)}-1.
  \ee
 Denote by $\mes G$ the Lebesgue measure of the set $G$. Since any half-interval
 $(y,y+\tau]$ contains at most $M(E[\tau/L]+1)$ points of the set $A$, we get
   $$
   \mes\left(\bigcup_{j=0}^N [(U_\e-j\tau)\cap(0,\,\tau)]\right)\le
 \sum_{j=0}^N \mes[U_\e\cap(j\tau,\,(j+1)\tau)]
   \le (N+1)2\e M(E[\tau/L]+1)<\tau.
  $$
 Hence there is $x\in(0,\,\tau)$ such that the points $x,\,x+\tau,\dots,x+N\tau$
 do not belong to $A_\e$.
 Therefore,
 $$
n(x+N\tau)-n(x)=\sum_{j=1}^N n(x+j\tau)-n(x+(j-1)\tau)=Np.
 $$
On the other hand, using (\ref{est}) with $h=N\tau$, we get
 $$
 |n(x+N\tau)-n(x)-N\tau|<4M.
 $$
Consequently, by (\ref{est-N}), we get $|\tau-p|<4M/N<16M^2\e/L$.

 Put
$\g(k)=a_k-k$ for all $k\in\Z$. We shall prove that
 \be\label{p}
 |\g(m+p)-\g(m)|<H\e\qquad\forall\,m\in\Z,
 \ee
with $H=5M+16M^2/L$. Suppose the contrary.  For example, let $\g(m+p)>\g(m)+H\e$ for
some $m\in\Z$. This yields that
 $$
 a_{m+p}>a_m+p+H\e>a_m+\tau+5M\e.
 $$
 Since $a_n\le a_m$ for $n<m$ and $a_n\ge a_{m+p}$ for $n>m+p$, we see that for all
 $t\in(a_m,\,a_m+5M\e)$ we have
 \be\label{num}
n(t+\tau)-n(t)= \card(A\cap(t,\,t+\tau])\le p-1.
 \ee
 On the other hand,
 $$
 \mes((a_m,\,a_m+L)\cap[U_\e\cup(U_\e-\tau)])\le2\e\,\card([a_m,\,a_m+L]\cap[A\cup(A-\tau)])\le
 4M\e<5M\e.
 $$
 Since the left-hand side of (\ref{num}) is equal to $p$ for all
 $t\in (a_m,\,a_m+L)\setminus[U_\e\cup(U_\e-\tau)]$,
we obtain a contradiction. In the same way we prove that the case $\g(m+p)<\g(m)-H\e$
is impossible as well. Hence (\ref{p}) is valid for all $m\in\Z$. If we continue the
function $\g$ as a linear function to each interval $(m,\,m+1)$, we obtain the
continuous function $f$ on $\R$ such that
 $$
 |f(x+p)-f(x)|<H\e,\quad\forall\,x\in\R.
 $$
 Since the number $p$ with this property exists in the $(16M^2\e/L)$-neighborhood of each
 $\e$-almost period $\tau$ of the  set $A$, we see that $f$ is an almost periodic
 function.  Theorem is proved. \ws

 \medskip
{\bf Signed multiple discrete sets.} Now we will consider some generalization of
discrete sets, namely signed multiple discrete sets in $\R^p$. This means that a
number $m(x)\in\Z\setminus\{0\}$ corresponds to each point $x$ from a discrete set. As
above, we denote this object by $A=\{(x,m(x))\}$ and the corresponding discrete set by
$s(A)$. Equality (\ref{f}) define the associate measure $\mu_A$ of $A$. Also, put
 $$
 A^+=\{(x,m(x)), x\in s(A),\,m(x)>0\},\qquad A^-=\{(x,m(x)), x\in s(A),\,m(x)<0\}.
 $$

In the case $p=2$ the definition coincides with the definition of the divisor of a
meromorphic function in the complex plane.

\medskip
 \begin{Def} A signed multiple discrete set $A$ is {\it almost periodic}, if its
  associate measure $\mu_A$ is almost periodic.
\end{Def}

Note that each continuous function with a compact support can be approximated by a
sequence of functions from $C^\infty$ with supports in a fixed ball. Therefore, if a
signed measure $\mu$ satisfies (\ref{b1}), we can take only functions $\f\in C^\infty$
in definition (\ref{conv}). Next, take a positive function $\f\in C^\infty$ such that
$\f(x)\equiv1$ for $|x|\le1$ and $\f(x)\equiv0$ for $|x|\ge2$ in (\ref{conv}). Since
almost periodic functions are bounded in $\R^p$, we see that every almost periodic in
the sense of distributions {\it positive} measure satisfies (\ref{b1}). Hence the
class of positive multiple discrete sets with almost periodic in the sense of
distributions associate measures coincides with the class of positive almost periodic
multiple discrete sets. But this assertion does not valid for signed multiple discrete
sets.

\begin{Th}\label{1}
  There is a signed multiple discrete set such that its associate measure is almost
  periodic in the sense of distributions and not almost periodic in the weak sense.
\end{Th}
{\bf Proof.} Let $\a(n)$, $n\in2\Z\setminus\{0\}$, be the greatest $k\in\N$ such that
$2^k$ is a divisor of $n$. Put
 $$
  a_n^+ =n+1/(\a(n)+1)^2,\quad a_n^-=n-1/(\a(n)+1)^2,\quad  n\in2\Z\setminus\{0\}.
 $$
  Define the signed multiple discrete set
$A=A^+\cup A^-$, where
 $$
     A^+=\{(a_n^+,\,\a(n))\}_{n\in2\Z\setminus\{0\}},\quad
     A^-=\{(a_n^-,\,\a(n))\}_{n\in2\Z\setminus\{0\}}.
 $$
The measure $\mu_A$ does not satisfy (\ref{b1}), therefore it is not almost periodic.
Let us show that $\mu_A$ is almost periodic in the sense of distributions.

Take a function $\varphi \in C^{\infty}$ such that
$\supp\,\varphi\subset(-1/2,\,1/2)$. Suppose that $\tau=2^pk$ for some
$p\in\N,\,k\in\Z$. If $|x-n|\ge3/4$ for all $n\in2\Z$, then the same is valid for the
point $x+\tau$, therefore,
 $$
 (\varphi*\mu_A)(x+\tau)=(\varphi*\mu_A)(x)=0.
 $$

If $|x-n|<3/4$ for $n\in\{2\Z:\,\a(n)\ge p\}\cup\{0\}$, then either
$(\varphi*\mu_A)(x)=0$, or
 $$
   |\varphi*\mu_A(x)|=\a(n)|\varphi(a_n^+ +x)
    -\varphi(a_n^-+x)|\le\a(n)M|a_n^+ -a_n^-|<2M/p,
 $$
where $M=\sup_\R|\varphi'(x)|$. Moreover, if this is the case, then also
$n+\tau\in\{2\Z:\,\a(n)\ge p\}\cup\{0\}$. Hence the same bound is valid for the value
$(\varphi*\mu_A)(x+\tau)$. We obtain
 $$
   |\varphi*\mu_A(x)-\varphi*\mu_A(x+\tau)|<4M/p.
 $$

Finally, if $|x-n|<3/4$ for $n\in\{2\Z:\,\a(n)<p\}$, then $\a(n+\tau)=\a(n)$.
Therefore, we have $a_{n+\tau}^\pm=a_n^\pm+\tau$, and
 $$
 (\varphi*\mu_A)(x+\tau)=\a(n+\tau)[\varphi(x+\tau-a_{n+\tau}^+)-\varphi(x+\tau-a_{n+\tau}^-)]
 $$
 $$
 =\a(n)[\varphi(x-a_n^+)-\varphi(x-a_n^-)]=(\varphi*\mu_A)(x).
 $$
 Consequently, all multiplies of $2^p$  are  $(4M/p)$-almost periods of the function
 $(\varphi*\mu_A)(x)$.  Thus $\mu_A$ is almost periodic in the sense  of distributions.
 Theorem is  proved.  \ws

\medskip
Note that we can check almost periodicity of function (\ref{conv}) only for positive
continuous functions $\f$ with an arbitrary small diameter of its support. Hence if
$A$ is a signed \A\ and
 $$
  \inf\{|x-y|:\,x\in s(A^+),\,y\in s(A^-)\}>0,
  $$
  then $A^+$ and $A^-$ are \A s as well. But this is false in the general case.
\begin{Th}\label{2}
  There is a signed almost periodic set $A$ such that $A^+$ and $A^-$ are not almost periodic.
\end{Th}
{\bf Proof.} Put
 $$
    A^+=\{(a_n^+, \,1)\}_{n\in2\Z},\quad  A^-=\{(a_n^-, \,-1) \}_{n\in2\Z},\quad A=A^+\cup A^-,
 $$
where points $a_n^\pm$ are the same as in the proof of Theorem \ref{1}.

 Following the proof of
Theorem \ref{1}, we can assure that $\mu_A$ is almost periodic in the sense of
distributions. Since the measure $\mu_A$ satisfies condition (\ref{b1}), we get that
$A$ is a signed \A.

We will use Theorem \ref{FK} for proving that $A^+$ is not almost periodic. Clearly,
the distance between any two points of $A^+$ has the form $2m+\b,\,m\in\N,\,|\b|<1/4$.
Hence whenever $\tau$ is an $\e$-almost period of $A^+$, $\e<1/4$, we have
$\tau=2n_0+\g,\,n_0\in\Z,\,|\g|<1/2$. But $0\not\in A^+$, hence the distance between
the point $2n_0+(\a(2n_0)+1)^{-2}-\tau$ and any point of $A^+$ is more than $1$. We
obtain a contradiction. Consequently, $A^+$ is not almost periodic. Analogously, $A^-$
is not almost periodic as well. Theorem is proved. \ws

Clearly, the measure $|\mu_A|=2\mu_{A^+}-\mu_A$ is not almost periodic as well.
Therefore we obtain
 \begin{Cor} The positive discrete set
$\{a_n^+\}\cup\{a_n^-\}$ does not almost periodic.
 \end{Cor}.

\end{document}